# On moment-density estimation in some biased models

## Robert M. Mnatsakanov[1] and Frits H. Ruymgaart[2]

*West Virginia University and Texas Tech University*

**Abstract:** This paper concerns estimating a probability density function $f$ based on iid observations from $g(x) = W^{-1} w(x) f(x)$, where the weight function $w$ and the total weight $W = \int w(x) f(x) \, dx$ may not be known. The length-biased and excess life distribution models are considered. The asymptotic normality and the rate of convergence in mean squared error (MSE) of the estimators are studied.

## 1. Introduction and preliminaries

It is known from the famous "moment problem" that under suitable conditions a probability distribution can be recovered from its moments. In Mnatsakanov and Ruymgaart [5, 6] an attempt has been made to exploit this idea and estimate a cdf or pdf, concentrated on the positive half-line, from its empirical moments.

The ensuing density estimators turned out to be of kernel type with a convolution kernel, provided that convolution is considered on the positive half-line with multiplication as a group operation (rather than addition on the entire real line). This does not seem to be unnatural when densities on the positive half-line are to be estimated; the present estimators have been shown to behave better in the right hand tail (at the level of constants) than the traditional estimator (Mnatsakanov and Ruymgaart [6]).

Apart from being an alternative to the usual density estimation techniques, the approach is particularly interesting in certain inverse problems, where the moments of the density of interest are related to those of the actually sampled density in a simple explicit manner. This occurs, for instance, in biased sampling models. In such models the pdf $f$ (or cdf $F$) of a positive random variable $X$ is of actual interest, but one observes a random sample $Y_1, \ldots, Y_n$ of $n$ copies of a random variable $Y$ with density

$$(1.1) \qquad g(y) = \frac{1}{W} \, w(y) \, f(y), \; y \geq 0,$$

where the weight function $w$ and the total weight

$$(1.2) \qquad W = \int_0^\infty w(x) \, f(x) \, dx,$$

[1]Department of Statistics, West Virginia University, Morgantown, WV 26506, USA, e-mail: rmnatsak@stat.wvu.edu
[2]Department of Mathematics & Statistics, Texas Tech University, Lubbock, TX 79409, USA, e-mail: h.ruymgaart@ttu.edu

may not be known. In this model one clearly has the relation

$$(1.3) \quad \mu_{k,F} = \int_0^\infty x^k f(x)\, dx = W \int_0^\infty y^k \frac{1}{w(y)}\, g(y)\, dy, \ k = 0, 1, \ldots,$$

and unbiased $\sqrt{n}$-consistent estimators of the moments of $F$ are given by

$$(1.4) \quad \widehat{\mu}_k = \frac{W}{n} \sum_{i=1}^n Y_i^k \frac{1}{w(Y_i)}.$$

If $w$ and $W$ are unknown they have to be replaced by estimators to yield $\widehat{\widehat{\mu}}_k$, say. In Mnatsakanov and Ruymgaart [7] moment-type estimators for the cdf $F$ of $X$ were constructed in biased models. In this paper we want to focus on estimating the density $f$ and related quantities. Following the construction pattern in Mnatsakanov and Ruymgaart [6], substitution of the empirical moments $\widehat{\mu}_k$ in the inversion formula for the density yields the estimators

$$(1.5) \quad \hat{f}_\alpha(x) = \frac{W}{n} \sum_{i=1}^n \frac{1}{w(Y_i)} \cdot \frac{\alpha - 1}{x \cdot (\alpha - 1)!} \cdot \left(\frac{\alpha}{x} Y_i\right)^{\alpha - 1} \exp(-\frac{\alpha}{x} Y_i), \ x \geq 0,$$

after some algebraic manipulation, where $\alpha$ is positive integer with $\alpha = \alpha(n) \to \infty$, as $n \to \infty$, at a rate to be specified later. If $W$ or $w$ are to be estimated, the empirical moments $\widehat{\widehat{\mu}}_k$ are substituted and we arrive at $\hat{\hat{f}}_\alpha$, say.

A special instance of model (1.1) to which this paper is devoted for the most part is length-biased sampling, where

$$(1.6) \quad w(y) = y, \ y \geq 0.$$

Bias and MSE for the estimator (1.5) in this particular case are considered in Section 3 and its asymptotic normality in Section 4. Although the weight function $w$ is known, its mean $W$ still remains to be estimated in most cases, and an estimator of $W$ is also briefly discussed. The literature on length-biased sampling is rather extensive; see, for instance Vardi [9], Bhattacharyya *et al.* [1] and Jones [4].

Another special case of (1.1) occurs in the study of the distribution of the excess of a renewal process; see, for instance, Ross [8] for a brief introduction. In this situation, it turns out that the sampled density satisfies (1.1) with

$$(1.7) \quad w(y) = \frac{1 - F(y)}{f(y)} = \frac{1}{h_F(y)}, \ y \geq 0,$$

where $h_F$ is the hazard rate of $F$. Although apparently $w$ and hence $W$ are not known here, they depend exclusively on $f$. In Section 5 we will briefly discuss some estimators for $f$, $h_F$ and $W$ and in particular show that they are all related to estimators of $g$ and its derivative. Estimating this $g$ is a "direct" problem and can formally be considered as a special case of (1.1) with $w(y) = 1$, $y \geq 0$ and $W = 1$. Investigating rates of convergence of the corresponding estimators is beyond the scope of this paper. Finally, in Section 6 we will compare the mean squared errors of the moment-density estimator $\hat{f}_\alpha^*$ introduced in the Section 2 and the kernel-density estimator $f_h$ studied by Jones [4] for the length-biased model. Throughout the paper let us denote by $G(a, b)$ a gamma distribution with shape and scale parameters $a$ and $b$, respectively. We carried out simulations for length-biased model (1.1) with $g$ as the gamma $G(2, 1/2)$ density and constructed corresponding graphs for $\hat{f}_\alpha^*$ and $f_h$. Also we compare the performance of the moment-type and kernel-type estimators for the model with excess life-time distribution when the target distribution $F$ is gamma $G(2, 2)$.

## 2. Construction of moment-density estimators and assumptions

Let us consider the general weighted model (1.1) and assume that the weight function $w$ is known. The estimated total weight $\hat{W}$ can be defined as follows:

$$\hat{W} = \Big(\frac{1}{n}\sum_{j=1}^{n} \frac{1}{w(Y_j)}\Big)^{-1}.$$

Substitution of the empirical moments

$$\hat{\hat{\mu}}_k = \frac{\hat{W}}{n}\sum_{i=1}^{n} Y_i^k \frac{1}{w(Y_i)}$$

in the inversion formula for the density (see, Mnatsakanov and Ruymgaart [6]) yields the construction

$$(2.1) \qquad \hat{\hat{f}}_\alpha(x) = \frac{\hat{W}}{n}\sum_{i=1}^{n} \frac{1}{w(Y_i)} \cdot \frac{\alpha - 1}{x \cdot (\alpha - 1)!} \cdot \Big(\frac{\alpha}{x}Y_i\Big)^{\alpha - 1} \exp(-\frac{\alpha}{x}Y_i).$$

Here $\alpha$ is positive integer and will be specified later. Note that the estimator $\hat{\hat{f}}_\alpha$ is the probability density itself. Note also that

$$\hat{W} = W + O_p(\frac{1}{\sqrt{n}}), \ n \to \infty,$$

(see, Cox [2] or Vardi [9]). Hence one can replace $\hat{W}$ in (2.1) by $W$.

Investigating the length-biased model, modify the estimator $\hat{\hat{f}}_\alpha$ and consider

$$\widehat{f}_\alpha^*(x) = \frac{W}{n}\sum_{i=1}^{n} \frac{1}{Y_i} \cdot \frac{\alpha}{x \cdot (\alpha - 1)!} \cdot \Big(\frac{\alpha}{x}Y_i\Big)^{\alpha - 1} \exp(-\frac{\alpha}{x}Y_i)$$

$$= \frac{1}{n}\sum_{i=1}^{n} \frac{W}{Y_i^2} \cdot \frac{1}{\Gamma(\alpha)} \cdot \Big(\frac{\alpha Y_i}{x}\Big)^{\alpha} \cdot \exp(-\frac{\alpha}{x}Y_i) = \frac{1}{n}\sum_{i=1}^{n} M_i.$$

In Sections 3 and 4 we will assume that the density $f$ satisfies

$$(2.2) \qquad f \in C^{(2)}([0,\infty)), \quad \text{with} \quad \sup_{t \geq 0}\big|f''(t)\big| = M < \infty.$$

In Section 5 we will consider the estimator of the unknown survival function $S$ given the sample $Y_1, \ldots, Y_n$ from pdf (1.1) where $w = (1 - F)/f$. Namely, we will use the moment-density estimator proposed in Mnatsakanov and Ruymgaart [6] which yields the estimator of $S = 1 - F$:

$$\hat{S}_\alpha(x) = \frac{W}{n}\sum_{i=1}^{n} \frac{1}{Y_i} \cdot \frac{1}{(\alpha - 1)!} \cdot \Big(\frac{\alpha}{x}Y_i\Big)^{\alpha} \exp(-\frac{\alpha}{x}Y_i).$$

We will assume that $F$ satisfies the conditions

$$(2.3) \qquad F \in C^{(2)}([0,\infty)), \quad \text{with} \quad \sup_{t \geq 0}\big|F''(t)\big| = L < \infty.$$

Throughout the paper the moment-type estimators will be considered at a fixed point $x > 0$, where $f(x) > 0$.

## 3. The bias and MSE of $\hat{f}^*_\alpha$

To study the asymptotic properties of $\widehat{f}^*_\alpha$ let us introduce for each $k \in \mathbb{N}$ the sequence of gamma $G(k(\alpha - 2) + 2, x/k\alpha)$ density functions

$$h_{\alpha,x,k}(u) = \frac{1}{\{k(\alpha - 2) + 1\}!} \left(\frac{k\alpha}{x}\right)^{k(\alpha-2)+2} u^{k(\alpha-2)+1} \tag{3.1}$$
$$\times \exp(-\frac{k\alpha}{x}u), \quad u \geq 0,$$

with mean $\{k(\alpha - 2) + 2\}x/(k\alpha)$ and variance $\{k(\alpha - 2) + 2\}x^2/(k\alpha)^2$. For each $k \in \mathbb{N}$, moreover, these densities form as well a delta sequence. Namely,

$$\int_0^\infty h_{\alpha,x,k}(u) f(u) \, du \to f(x), \text{ as } \alpha \to \infty,$$

uniformly on any bounded interval (see, for example, Feller [3], vol. II, Chapter VII). This property of $h_{\alpha,x,k}$, when $k = 2$ is used in (3.10) below. In addition, for $k = 1$ we have

$$\int_0^\infty u \, h_{\alpha,x,1}(u) \, du = x, \tag{3.2}$$

$$\int_0^\infty (u - x)^2 h_{\alpha,x,1}(u) \, du = \frac{x^2}{\alpha}. \tag{3.3}$$

**Theorem 3.1.** *Under the assumptions (2.2) the bias of $\widehat{f}^*_\alpha$ satisfies*

$$E\widehat{f}^*_\alpha(x) - f(x) = \frac{x^2 f''(x)}{2 \cdot \alpha} + o\left(\frac{1}{\alpha}\right), \quad \text{as } \alpha \to \infty. \tag{3.4}$$

*For the Mean Squared Error* (MSE) *we have*

$$MSE\{\widehat{f}^*_\alpha(x)\} = n^{-4/5} \left[\frac{W \cdot f(x)}{2\sqrt{\pi}x^2} + \frac{x^4\{f''(x)\}^2}{4}\right] + o(1), \tag{3.5}$$

*provided that we choose* $\alpha = \alpha(n) \sim n^{2/5}$.

*Proof.* Let $M_i = W \cdot Y_i^{-1} \cdot h_{\alpha,x,1}(Y_i)$. Then

$$E M_i^k = \int_0^\infty W^k \cdot Y_i^{-k} h^k_{\alpha,x,1}(y) g(y) \, dy$$
$$= \int_0^\infty \frac{W^k}{\{y \cdot (\alpha - 1)!\}^k} \left(\frac{\alpha}{x}\right)^{k\alpha} y^{k(\alpha-1)} \exp\left(-\frac{k\alpha}{x}y\right) \frac{y \cdot f(y)}{W} dy \tag{3.6}$$
$$= W^{k-1} \int_0^\infty \frac{1}{\{(\alpha - 1)!\}^k} \left(\frac{\alpha}{x}\right)^{k\alpha} y^{k(\alpha-2)+1} \exp\left(-\frac{k\alpha}{x}y\right) f(y) \, dy$$
$$= W^{k-1} \left(\frac{\alpha}{x}\right)^{2(k-1)} \frac{\{k(\alpha - 2) + 1\}!}{\{(\alpha - 1)!\}^k} \frac{1}{k^{k(\alpha-2)+2}} \int_0^\infty h_{\alpha,x,k}(y) f(y) dy.$$

In particular, for $k = 1$:

$$E\widehat{f}^*_\alpha(x) = f_\alpha(x) = W \int_0^\infty \frac{1}{y^2} \cdot \frac{1}{\Gamma(\alpha)} \cdot \left(\frac{\alpha}{x}y\right)^\alpha \exp(-\frac{\alpha}{x}y) \frac{yf(y)}{W} dy \tag{3.7}$$
$$= \int_0^\infty h_{\alpha,x,1}(y) f(y) dy = E M_i.$$

This yields for the bias ($\mu = x, \sigma^2 = x^2/\alpha$)

$$
\begin{aligned}
f_\alpha(x) - f(x) &= \int_0^\infty h_{\alpha,x,1}(y)\{f(y) - f(x)\}du \\
&= \int_0^\infty h_{\alpha,x,1}(y)\{f(x) + (y-x)f'(x) \\
&\quad + \frac{1}{2}\int_0^\infty (y-x)^2\{f''(\tilde{y}) - f(x)\}dy \\
&= \frac{1}{2}\int_0^\infty (y-x)^2 h_{\alpha,x,1}(y) f''(x) du \\
&\quad + \frac{1}{2}\int_0^\infty (y-x)^2 h_{\alpha,x,1}(y)\{f''(\tilde{y}) - f''(x)\}dy \\
&= \frac{1}{2}\frac{x^2}{\alpha}f''(x) + o\left(\frac{1}{\alpha}\right), \quad \text{as } \alpha \to \infty.
\end{aligned}
\tag{3.8}
$$

For the variance we have

$$
\operatorname{Var}\widehat{f}_\alpha^*(x) = \frac{1}{n}\operatorname{Var} M_i = \frac{1}{n}\{E M_i^2 - f_\alpha^2(x)\}. \tag{3.9}
$$

Applying (3.6) for $k = 2$ yields

$$
\begin{aligned}
E M_i^2 &= W\frac{\alpha^2}{x^2}\frac{(2\alpha-3)!}{\{(\alpha-1)!\}^2}\frac{1}{2^{2\alpha-2}}\int_0^\infty h_{\alpha,x,2}(u)f(u)du \\
&\sim \frac{\alpha^2}{x^2}\frac{W}{\sqrt{2\pi}}\frac{e^{-(2\alpha-3)}\{(2\alpha-3)\}^{(2\alpha-3)+1/2}}{e^{-2(\alpha-1)}\{(\alpha-1)\}^{2(\alpha-1)+1}}\frac{1}{2^{2(\alpha-1)}} \\
&\quad \times \int_0^\infty h_{\alpha,x,2}(u)f(u)du = \frac{W}{2\sqrt{\pi}}\frac{\sqrt{\alpha}}{x^2}\int_0^\infty h_{\alpha,x,2}(u)f(u)du \\
&= \frac{W}{2\sqrt{\pi}}\frac{\sqrt{\alpha}}{x^2}\{f(x) + o(1)\} = \frac{W}{2\sqrt{\pi}}\frac{\sqrt{\alpha}}{x^2}f(x) + o(\sqrt{\alpha})
\end{aligned}
\tag{3.10}
$$

as $\alpha \to \infty$. Now inserting this in (3.9) we obtain

$$
\begin{aligned}
\operatorname{Var}\widehat{f}_\alpha^*(x) &= \frac{1}{n}\left[\frac{W}{2\sqrt{\pi}}\frac{\sqrt{\alpha}}{x^2}f(x) + o(\sqrt{\alpha}) - \left\{f(x) + O\left(\frac{1}{\alpha}\right)\right\}^2\right] \\
&= \frac{W\sqrt{\alpha}}{2n\sqrt{\pi}}\frac{f(x)}{x^2} + o\left(\frac{\sqrt{\alpha}}{n}\right).
\end{aligned}
\tag{3.11}
$$

Finally, this leads to the MSE of $\widehat{f}_\alpha^*(x)$:

$$
MSE\{\widehat{f}_\alpha^*(x)\} = \frac{W\sqrt{\alpha}}{2n\sqrt{\pi}}\frac{f(x)}{x^2} + \frac{1}{4}\frac{x^4}{\alpha^2}\{f''(x)\}^2 + o\left(\frac{\sqrt{\alpha}}{n}\right) + o\left(\frac{1}{\alpha^2}\right). \tag{3.12}
$$

For optimal rate we may take

$$
\alpha = \alpha_n \sim n^{2/5}, \tag{3.13}
$$

assuming that $n$ is such that $\alpha_n$ is an integer. By substitution (3.13) in (3.12) we find (3.5). □

**Corollary 3.1.** *Assume that the parameter* $\alpha = \alpha(x)$ *is chosen locally for each* $x > 0$ *as follows*

$$(3.14) \qquad \alpha(x) = n^{2/5} \cdot \{\frac{\pi}{4 \cdot W^2}\}^{1/5} \left[\frac{x^3 \cdot f''(x)}{\sqrt{f(x)}}\right]^{4/5}, \; f''(x) \neq 0.$$

*Then the estimator* $\widetilde{f}^*_\alpha(x) = \widehat{f}^*_\alpha(x)$ *satisfies*

$$(3.15) \qquad \widetilde{f}^*_\alpha(x) = \frac{1}{n} \sum_{i=1}^n \frac{W}{Y_i^2} \cdot \frac{1}{\Gamma(\alpha(x))} \cdot \left(\frac{\alpha(x)}{x} Y_i\right)^{\alpha(x)} \exp\{-\frac{\alpha(x)}{x} Y_i\}$$

$$(3.16) \quad MSE\{\widetilde{f}^*_\alpha(x)\} = n^{-4/5} \left[\frac{W^2 \cdot f''(x) \cdot f^2(x)}{\pi \cdot x^2 \sqrt{2}}\right]^{2/5} + o(1), \text{ as } n \to \infty.$$

*Proof.* Assuming the first two terms in the right hand side of (3.12) are equal to each other one obtains that for each $n$ the function $\alpha = \alpha(x)$ can be chosen according to (3.14). This yields the proof of Corollary 1. □

## 4. The asymptotic normality of $\widehat{f}^*_\alpha$

Now let us derive the limiting distributions of $\widehat{f}^*_\alpha$. The following statement is valid.

**Theorem 4.1.** *Under the assumptions* (2.2) *and* $\alpha = \alpha(n) \sim n^\delta$, *for any* $0 < \delta < 2$, *we have, as* $\alpha \to \infty$,

$$(4.1) \qquad \frac{\widehat{f}^*_\alpha(x) - f_\alpha(x)}{\sqrt{Var\widehat{f}^*_\alpha(x)}} \to_d Normal(0,1).$$

*Proof.* Let $0 < C < \infty$ denote a generic constant that does not depend on $n$ but whose value may vary from line to line. Note that for arbitrary $k \in \mathbb{N}$ the "$c_r$-inequality" entails that $E|M_i - f_\alpha(x)|^k \leq C \, E\, M_i^k$, in view of (3.6) and (3.7). Now let us choose the integer $k > 2$. Then it follows from (3.6) and (3.11) that

$$(4.2) \qquad \frac{\sum_{i=1}^n E|\frac{1}{n}\{M_i - f_\alpha(x)\}|^k}{\{Var\widehat{f}_\alpha(x)\}^{k/2}} \leq C \frac{n^{1-k} k^{-1/2} \alpha^{k/2 - 1/2}}{(n^{-1} \alpha^{1/2})^{k/2}}$$

$$= C \frac{1}{\sqrt{k}} \frac{\alpha^{k/4 - 1/2}}{n^{k/2 - 1}} \to 0, \text{ as } n \to \infty,$$

for $\alpha \sim n^\delta$. Thus the Lyapunov's condition for the central limit theorem is fulfilled and (4.1) follows for any $0 < \delta < 2$. □

**Theorem 4.2.** *Under the assumptions* (2.2) *we have*

$$(4.3) \qquad \frac{n^{1/2}}{\alpha^{1/4}} \{\widehat{f}^*_\alpha(x) - f(x)\} \to_d Normal\left(0, \frac{W \cdot f(x)}{2 \, x^2 \sqrt{\pi}}\right),$$

*as* $n \to \infty$, *provided that we take* $\alpha = \alpha(n) \sim n^\delta$ *for any* $\frac{2}{5} < \delta < 2$.

*Proof.* This is immediate from (3.11) and (4.1), since combined with (3.8) entails that $n^{1/2} \alpha^{-1/4} \{f_\alpha(x) - f(x)\} = O(n^{1/2} \alpha^{-5/4}) = O(n^{-\frac{5\delta - 2}{4}}) = o(1)$, as $n \to \infty$, for the present choice of $\alpha$. □

**Corollary 4.1.** *Let us assume that* (2.2) *is valid. Consider* $\widetilde{f}_\alpha^*(x)$ *defined in* (3.15) *with* $\alpha(x)$ *given by* (3.14). *Then*

$$(4.4) \qquad \frac{n^{1/2}}{\alpha(x)^{1/4}} \{\widetilde{f}_\alpha^*(x) - f(x)\} \to_d \text{Normal}\left([\frac{W f(x)}{2 x^2 \sqrt{\pi}}]^{1/2}, \frac{W f(x)}{2 x^2 \sqrt{\pi}}\right),$$

*as* $n \to \infty$ *and* $f''(x) \neq 0$.

*Proof.* From (4.1) and (3.11) with $\alpha = \alpha(x)$ defined in (3.13) it is easy to see that

$$(4.5) \qquad \frac{n^{1/2}}{\alpha(x)^{1/4}} \{\widetilde{f}_\alpha^*(x) - E\widetilde{f}_\alpha^*(x)\} = \text{Normal}\left(0, \frac{W f(x)}{2 x^2 \sqrt{\pi}}\right) + o_P(\frac{1}{n^{2/5}}),$$

as $n \to \infty$. Application of (3.4) where $\alpha = \alpha(x)$ is defined by (3.14) yields (4.4). □

**Corollary 4.2.** *Let us assume that* (2.2) *is valid. Consider* $\widetilde{f}_{\alpha^*}^*(x)$ *defined in* (3.15) *with* $\alpha^*(x)$ *given by*

$$(4.6) \qquad \alpha^*(x) = n^\delta \cdot \{\frac{\pi}{4 \cdot W^2}\}^{1/5} \left[\frac{x^3 \cdot f''(x)}{\sqrt{f(x)}}\right]^{4/5}, \quad \frac{2}{5} < \delta < 2.$$

*Then when* $f''(x) \neq 0$, *and letting* $n \to \infty$, *it follows that*

$$(4.7) \qquad \frac{n^{1/2}}{\alpha^*(x)^{1/4}} \{\widetilde{f}_{\alpha^*}^*(x) - f(x)\} \to_d \text{Normal}\left(0, \frac{W f(x)}{2 x^2 \sqrt{\pi}}\right).$$

*Proof.* Again from (4.1) and (3.11) with $\alpha = \alpha^*(x)$ defined in (4.6) it is easy to see that

$$(4.8) \qquad \frac{n^{1/2}}{\alpha^*(x)^{1/4}} \{\widetilde{f}_{\alpha^*}^*(x) - E\widetilde{f}_{\alpha^*}^*(x)\} = \text{Normal}\left(0, \frac{W f(x)}{2 x^2 \sqrt{\pi}}\right) + o_P(1),$$

as $n \to \infty$. On the other hand application of (3.4) where $\alpha = \alpha^*(x)$ is defined by (4.6) yields

$$(4.9) \qquad \frac{n^{1/2}}{\alpha^*(x)^{1/4}} \{E\widetilde{f}_{\alpha^*}^*(x) - f(x)\} = O\left(\frac{C(x)}{n^{(5\delta-2)/4}}\right),$$

as $n \to \infty$. Here $C(x) = \{\frac{W f(x)}{2 x^2 \sqrt{\pi}}\}^{1/2}$. Combining (4.8) and (4.9) yields (4.7). □

## 5. An application to the excess life distribution

Assume that the random variable $X$ has cdf $F$ and pdf $f$ defined on $[0, \infty)$ with $F(0) = 0$. Denote the hazard rate function $h_F = f/S$, where $S = 1 - F$ is the corresponding survival function of $X$. Assume also that the sampled density $g$ satisties (1.1) and (1.7). It follows that

$$(5.1) \qquad g(y) = \frac{1}{W} \{1 - F(y)\}, \quad y \geq 0.$$

It is also immediate that $W = 1/g(0)$ and, $f(y) = -W g'(y) = -\frac{g'(y)}{g(0)}$, $y \geq 0$, so that

$$h_F(y) = -\frac{g'(y)}{g(y)}, \quad y \geq 0.$$

Suppose now that we are given $n$ independent copies $Y_1, \ldots, Y_n$ of a random variable $Y$ with cdf $G$ and density $g$ from (5.1).

To recover $F$ or $S$ from the sample $Y_1, \ldots, Y_n$ use the moment-density estimator from Mnatsakanov and Ruymgaart [6], namely

$$(5.2) \qquad \hat{\hat{S}}_\alpha(x) = \frac{\hat{W}}{n} \sum_{i=1}^n \frac{1}{Y_i} \cdot \frac{1}{(\alpha-1)!} \cdot \left(\frac{\alpha}{x} Y_i\right)^\alpha \exp\left(-\frac{\alpha}{x} Y_i\right).$$

Where the estimator $\hat{W}$ can be defined as follows:

$$\hat{W} = \frac{1}{\hat{g}(0)}.$$

Here $\hat{g}$ is any estimator of $g$ based on the sample $Y_1, \ldots, Y_n$.

**Remark 5.1.** As has been noted at the end of Section 1, estimating $g$ from $Y_1, \ldots, Y_n$ is a "direct" problem and an estimator of $g$ can be constructed from (1.5) with $W$ and $w(Y_i)$ both replaced by 1. This yields

$$(5.3) \qquad \hat{g}_\alpha(y) = \frac{1}{n} \sum_{i=1}^n \frac{1}{y} \cdot \frac{\alpha-1}{(\alpha-1)!} \cdot \left(\frac{\alpha}{y} Y_i\right)^{\alpha-1} \exp\left(-\frac{\alpha}{y} Y_i\right), \quad y \geq 0.$$

The relations above suggest the estimators

$$\hat{f}(y) = -\frac{\hat{g}'_\alpha(y)}{\hat{g}_\alpha(0)}, \quad y \geq 0,$$

$$\hat{h}_F(y) = -\frac{\hat{g}'_\alpha(y)}{\hat{g}_\alpha(y)}, \quad \hat{w}(y) = -\frac{\hat{g}_\alpha(y)}{\hat{g}'_\alpha(y)}, \quad y \geq 0.$$

Here let us assume for simplicity that $W$ is known and construct the estimator of survival function $S$ as follows:

$$(5.4) \qquad \hat{S}_\alpha(x) = \frac{1}{n} \sum_{i=1}^n \frac{W}{Y_i} \cdot \frac{1}{\Gamma(\alpha)} \cdot \left(\frac{\alpha}{x} Y_i\right)^\alpha \exp\left(-\frac{\alpha}{x} Y_i\right) = \frac{1}{n} \sum_{i=1}^n L_i.$$

**Theorem 5.1.** *Under the assumptions* (2.3) *the bias of $\hat{S}_\alpha$ satisfies*

$$(5.5) \qquad E\hat{S}_\alpha(x) - S(x) = -\frac{x^2 f'(x)}{2 \cdot \alpha} + o\left(\frac{1}{\alpha}\right), \quad \text{as } \alpha \to \infty.$$

*For the Mean Squared Error* (MSE) *we have*

$$(5.6) \qquad MSE\{\hat{S}_\alpha(x)\} = n^{-4/5} \left[\frac{W \cdot S(x)}{2 \cdot x\sqrt{\pi}} + \frac{x^4 \{f'(x)\}^2}{4}\right] + o(1),$$

*provided that we choose $\alpha = \alpha(n) \sim n^{2/5}$.*

*Proof.* By a similar argument to the one used in (3.8) and (3.10) it can be shown that

$$(5.7) \qquad E\hat{S}_\alpha(x) - S(x) = \int_0^\infty h_{\alpha,x,1}(u)\{S(u) - S(x)\}du$$

$$= -\frac{1}{2}\frac{x^2}{\alpha}f'(x) + o\left(\frac{1}{\alpha}\right), \quad \text{as } \alpha \to \infty$$

and

$$(5.8) \quad E L_i^2 = \frac{W}{2\sqrt{\pi}} \frac{\sqrt{\alpha}}{x} S(x) + o(\sqrt{\alpha}), \text{ as } \alpha \to \infty,$$

respectively. So that combining (5.7) and (5.8) yields (5.6). □

**Corollary 5.1.** *If the parameter $\alpha = \alpha(x)$ is chosen locally for each $x > 0$ as follows*

$$(5.9) \quad \alpha(x) = n^{2/5} \cdot \{\frac{\pi}{4 \cdot W^2}\}^{1/5} \cdot x^2 \cdot \left[\frac{f'(x)}{\sqrt{1 - F(x)}}\right]^{4/5}, \, f'(x) \neq 0.$$

*then the estimator (5.4) with $\alpha = \alpha(x)$ satisfies*

$$MSE\{\hat{S}_\alpha(x)\} = n^{-4/5} \left[\frac{W^2 \cdot f'(x) \cdot (1 - F(x))^2}{\pi \sqrt{2}}\right]^{2/5} + o(1), \text{ as } n \to \infty.$$

**Theorem 5.2.** *Under the assumptions (2.3) and $\alpha = \alpha(n) \sim n^\delta$ for any $0 < \delta < 2$ we have, as $n \to \infty$,*

$$(5.10) \quad \frac{\hat{S}_\alpha(x) - E\hat{S}_\alpha(x)}{\sqrt{Var\, \hat{S}_\alpha(x)}} \to_d Normal(0, 1).$$

**Theorem 5.3.** *Under the assumptions (2.3) we have*

$$(5.11) \quad \frac{n^{1/2}}{\alpha^{1/4}} \{\hat{S}_\alpha(x) - S(x)\} \to_d Normal\left(0, \frac{W \cdot S(x)}{2x\sqrt{\pi}}\right),$$

*as $n \to \infty$, provided that we take $\alpha = \alpha(n) \sim n^\delta$ for any $\frac{2}{5} < \delta < 2$.*

**Corollary 5.2.** *If the parameter $\alpha = \alpha(x)$ is chosen locally for each $x > 0$ according to (5.9) then for $\hat{S}_\alpha(x)$ defined in (5.4) we have*

$$\frac{n^{1/2}}{\alpha(x)^{1/4}} \{\hat{S}_\alpha(x) - S(x)\} \to_d Normal\left(-[\frac{W\, S(x)}{2x\sqrt{\pi}}]^{1/2}, \frac{W\, S(x)}{2x\sqrt{\pi}}\right),$$

*provided $f'(x) \neq 0$ and $n \to \infty$.*

**Corollary 5.3.** *If the parameter $\alpha = \alpha^*(x)$ is chosen locally for each $x > 0$ according to*

$$(5.12) \quad \alpha^*(x) = n^\delta \cdot \{\frac{\pi}{4 \cdot W^2}\}^{1/5} \cdot x^2 \cdot \left[\frac{f'(x)}{\sqrt{1 - F(x)}}\right]^{4/5}, \, \frac{2}{5} < \delta < 2,$$

*then for $\hat{S}_{\alpha^*}(x)$ defined in (5.4) we have*

$$\frac{n^{1/2}}{\alpha^*(x)^{1/4}} \{\hat{S}_{\alpha^*}(x) - S(x)\} \to_d Normal\left(0, \frac{W\, S(x)}{2x\sqrt{\pi}}\right),$$

*provided $f'(x) \neq 0$ and $n \to \infty$.*

Note that the proofs of all statements from Theorems 5.2 and 5.3 are similar to the ones from Theorems 4.1 and 4.2, respectively.

## 6. Simulations

At first let us compare the graphs of our estimator $\widehat{f}_\alpha^*$ and the kernel-density estimator $f_h$ proposed by Jones [4] in the length-biased model:

$$(6.1) \qquad f_h(x) = \frac{\hat{W}}{nh} \sum_{i=1}^n \frac{1}{Y_i} \cdot K\left(\frac{x-Y_i}{h}\right), \ x > 0.$$

Assume, for example, that the kernel $K(x)$ is a standard normal density, while the bandwidth $h = O(n^{-\beta})$, with $0 < \beta < 1/4$. Here $\hat{W}$ is defined as follows

$$\hat{W} = \left(\frac{1}{n}\sum_{j=1}^n \frac{1}{Y_j}\right)^{-1}.$$

In Jones [4] under the assumption that $f$ has two continuous derivatives, it was shown that as $n \to \infty$

$$(6.2)\ MSE\{f_h(x)\} = \operatorname{Var} f_h(x) + \operatorname{bias}^2\{f_h\}(x)$$
$$\sim \frac{Wf(x)}{nhx}\int_0^\infty K^2(u)du + \frac{1}{4}h^4\{f''(x)\}^2\left[\int_0^\infty u^2 K(u)du\right]^2.$$

Comparing (6.2) with (3.12), where $\alpha = h^{-2}$, one can see that the variance term $\operatorname{Var} \widehat{f}_\alpha^*(x)$ for the moment-density estimator could be smaller for large values of $x$ than the corresponding $\operatorname{Var}\{f_h(x)\}$ for the kernel-density estimator. Near the origin the variability of $f_h$ could be smaller than that of $\widehat{f}_\alpha^*$. The bias term of $\widehat{f}_\alpha^*$ contains the extra factor $x^2$, but as the simulations suggest this difference is compensated by the small variability of the moment-density estimator.

We simulated $n = 300$ copies of length-biased r.v.'s from gamma $G(2, 1/2)$. The corresponding curves for $f$ (solid line) and its estimators $\widehat{f}_\alpha^*$ (dashed line), and $f_h$ (dotted line), respectively are plotted in Figure 1. Here we chose $\alpha = n^{2/5}$ and

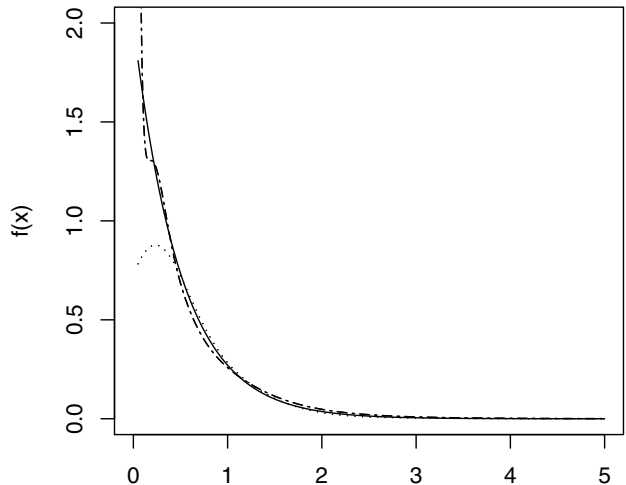

Figure 1.

FIG 1.

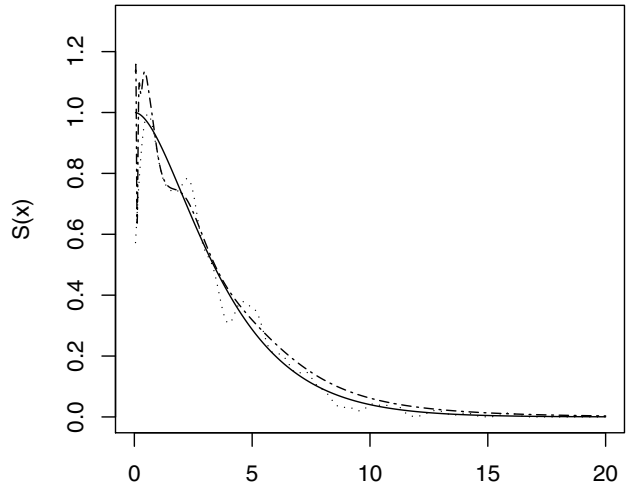

Figure 2.

Fig 2.

$h = n^{-1/5}$, respectively. To construct the graphs for the moment-type estimator $\hat{S}_\alpha$ defined by (5.4) and the kernel-type estimator $S_h$ defined in a similar way as the one given by (6.1) let us generate $n = 400$ copies of r.v.'s $Y_1, \ldots, Y_n$ with pdf $g$ from (5.1) with $W = 4$ and

$$1 - F(x) = e^{-\frac{x}{2}} + \frac{x}{2} e^{-\frac{x}{2}}, x \geq 0.$$

We generated $Y_1, \ldots, Y_n$ as a mixture of two gamma $G(1,2)$ and $G(2,2)$ distributions with equal proportions. In the Figure 2 the solid line represents the graph of $S = 1 - F$ while the dashed and dotted lines correspond to $\hat{S}_\alpha$ and $S_h$, respectively. Here again we have $\alpha = n^{2/5}$ and $h = n^{-1/5}$.